\documentclass[12pt,twoside]{article}

\usepackage{a4}
\usepackage{exscale}
\usepackage{fancyhdr}
\usepackage[centertags]{amsmath}
\usepackage{amsfonts}
\usepackage{amssymb}
\usepackage{theorem}
\usepackage{ulem}
\usepackage{isolatin1}
\usepackage{epsfig}

\renewcommand{\sectionmark}[1]
                    {\markboth{Kapitel \thesection\ #1}{}}
\renewcommand{\sectionmark}[1]
                 {\markright{} }

\newtheorem{thm}{Theorem}
\newtheorem{definition}{Definition}
\newtheorem{lem}[thm]{Lemma}
\newtheorem{corollary}[thm]{Corollary}

\newtheorem{proposition}[thm]{Proposition}

\numberwithin{equation}{section}

\def\parmod{\parskip=2pt plus1pt minus1pt}

\newenvironment{einr}{\parmod
                      \begin{list}{}
                        {\setlength{\rightmargin}{0cm}
                         \setlength{\leftmargin}{0,75cm}
                         \setlength{\labelwidth}{0cm}
                         \setlength{\parsep}{1pt}
                         \setlength{\itemsep}{1pt}
                         \setlength{\topsep}{1pt}
                         \setlength{\partopsep}{0pt}
                         \setlength{\labelsep}{0cm}
                         \setlength{\listparindent}{0pt}
                         \setlength{\itemindent}{0pt}}
                      \item[] \ignorespaces}
                     {\unskip \end{list}}

\newenvironment{proof}
        {\pagebreak[2] \vspace{-1pt}{\bf Proof.}  }
        {\hfill $\blacksquare$ \vspace{2pt}}

\newenvironment{proofof}[1]
        {\pagebreak[2] \vspace{-1pt}{\bf Proof of #1.}  }
        {\hfill $\blacksquare$ \vspace{2pt}}

\def\fa{\mathcal{F}}
\def\ha{\mathcal{H}}
\def\ma{\mathcal{M}}

\def\pa{\mathcal{P}}
\def\za{\mathcal{Z}}

\def\nat{{\rm I\! N}}

\def\co{{\mathbb C}}

\def\dk{{\mathbb D}}

\def\rd{\partial\dk}

\def\l{\left}
\def\r{\right}
\def\gl{\left\{}
\def\gr{\right\}}
\def\kl{\left(}
\def\kr{\right)}

\def\limn{\lim_{n\to\infty}}

\def\ti{\widetilde}

\def\In{\subseteq}

\def\mi{\setminus}
\def\abb{\longrightarrow}

\def\plog{\log^+}

\renewcommand{\rho}{\varrho}
\renewcommand{\phi}{\varphi}
\renewcommand{\epsilon}{\varepsilon}

\def\beq{\begin{equation}}
\def\eeq{\end{equation}}
\def\beqar{\begin{eqnarray}}
\def\eeqar{\end{eqnarray}}
\def\beqaro{\begin{eqnarray*}}
\def\eeqaro{\end{eqnarray*}}
\def\bsat{\begin{thm}}
\def\esat{\end{thm}}
\def\blem{\begin{lem}}
\def\elem{\end{lem}}
\def\bkor{\begin{corollary}}
\def\ekor{\end{corollary}}
\def\bprop{\begin{proposition}}
\def\eprop{\end{proposition}}
\def\bdefin{\begin{definition}}
\def\edefin{\end{definition}}
\def\bbew{\begin{proof}}
\def\ebew{\end{proof}}
\def\bbewo{\begin{proofof}}
\def\ebewo{\end{proofof}}
\def\beinr{\begin{einr}}
\def\eeinr{\end{einr}}

\renewcommand{\rho}{\varrho}
\renewcommand{\phi}{\varphi}

\begin{document}

\thispagestyle{plain}

\fancyhead[CE]{J\"urgen Grahl and Shahar Nevo}
\fancyhead[CO]{Spherical derivatives and normal families}
\fancyhead[LE,RO]{\thepage}
\fancyhead[LO,RE]{}
\fancyfoot[CE,CO]{}

\begin{center}
{\LARGE\bf A note on spherical derivatives and normal families \\[18pt]}

{\Large \it J\"urgen Grahl and Shahar Nevo}
\end{center}

\centerline{\bf Abstract}

We show that a family of meromorphic functions in the unit disk $\dk$
whose spherical derivatives are uniformly bounded away from zero is
normal. Furthermore, we show that for each $f$ meromorphic in $\dk$ we
have
$$\inf_{z\in\dk} f^\#(z)\le \frac{1}{2}$$
where $f^\#$ denotes the spherical derivative of $f$.

{\it 2000 Mathematics Subject Classification: 30A10, 30D45}

\section{Introduction and statement of results}

By a well-known result of Marty, a family $\fa$ of meromorphic
functions in a domain $D\In\co$ is normal (in the sense of
Montel) if and only if the family $\fa^\#:=\gl f^\#\;:\;f\in\fa\gr$ of
the corresponding spherical derivatives is locally uniformly bounded
in $D$; here, $f^\#$ is defined by $f^\#:=\frac{|f'|}{1+|f|^2}$.

To our best knowledge, it hasn't been studied so far what can be said if
$\fa^\#$ is uniformly bounded away from zero. It is the aim of the
present paper to tackle this question.

If $D$ is an arbitrary domain in $\co$, by $\ma(D)$ we denote
the space of all functions which are meromorphic in $D$ and by
$\ha(D)$ the space of all functions which are analytic in $D$. We
write $P_f:=f^{-1}(\gl\infty\gr)$ for the set of poles of a
meromorphic function $f$. If $\fa$ is some family of meromorphic
functions in a fixed domain, we
set
$$\fa':=\gl f'\;:\; f\in\fa\gr$$
for the family of the corresponding derivatives. Furthermore, we
denote the open (resp. closed) disk with center $c$ and radius $r$
by $U_r(c)$ (resp. $B_r(c)$) and set $\dk:=U_1(0)$ for the open unit
disk. Since normality is a local property, we can restrict all our
considerations concerning normal families to the unit disk.

Our main result is the following.

\bsat\label{normalresult}
Let some $\varepsilon>0$ be given and set
$$\fa:=\gl f\in\ma(\dk)\;:\; f^\#(z)\ge \varepsilon \mbox{ for all }
z\in\dk\gr.$$
Then $\fa$ is normal in $\dk$.
\esat

An immediate consequence is the existence of some new universal
constant for meromorphic functions.

\bkor \label{universal}
There exists some universal constant $C<\infty$ such that
$$\inf_{z\in\dk} f^\#(z)\le C$$
for each function $f$ meromorphic in $\dk$.
\ekor

\bbew{}
Otherwise, we could find a sequence $(f_n)_n$ in $\ma(\dk)$ such that
$$\inf_{z\in\dk} f_n^\#(z)\ge n$$
for each $n$. By Theorem \ref{normalresult}, $(f_n)_n$ is normal in
$\dk$. So by Marty's theorem there exists some constant $M<\infty$
such that $f_n^\#(z)\le M$ for all $n$ and all $z\in U_{1/2}(0)$. But
then of course $\inf_{z\in\dk} f_n^\#(z)\le M$ for all $n$. For $n>M$
this gives a contradiction.
\ebew

It's natural to conjecture that the extremal function for this problem
is the identity function and that it is more or less unique (up to
rigid motions of the Riemann sphere\footnote{It is well known that
  those rigid motions can be described as M\"obius transformations of
  the form $T(z)=\frac{az+b}{cz+d}$ such that the matrix
$\begin{pmatrix} a & b \\ c & d \end{pmatrix}$
is unitary.}). This conjecture turns out to be true.

\bsat \label{universalconstant}
Let $f$ be meromorphic in $\dk$. Then
\beq\label{infSphericalDeriv}
\inf_{z\in\dk} f^\#(z)\le \frac{1}{2}.
\eeq
Here, equality holds if and only if $f$ is a rigid motion of the
Riemann sphere.
\esat







In this context, it's interesting to note that the infimum of the
spherical derivative can only in trivial cases be attained in the
interior of the domain of definition, i.e. that there is a minimum
principle for the spherical derivative.

\bprop[Minimum principle for the spherical derivative]
\label{minimumprinciple}
Let $f$ be meromorphic in $\dk$. If $f^\#$ attains a local mimimum in
$z_0\in\dk$, then $f^\#(z_0)=0$.
\eprop

Since we couldn't find a reference for this fact, we give a short
proof (in Section \ref{proofs}) for the convenience of the reader.

An immediate consequence of Theorem \ref{normalresult}, Marty's
theorem, Proposition~\ref{minimumprinciple} and of some standard
compactness and continuity arguments is the following normality
criterion for exceptional functions of the spherical derivative.

\bkor{}
Let $h:\dk\abb [0;\infty)$ be a continuous and non-negative function
such that $\liminf_{z\to \zeta} h(z)>0$ for each $\zeta\in\rd$. Then
the family
$$\fa:=\gl f\in\ma(\dk)\;:\; f^\#(z)\ne h(z) \mbox{ for all }
z\in\dk\gr$$
is normal.
\ekor


\section{The Analytic Case of Theorem \ref{normalresult}} \label{analyticcase}

Throughout the paper, one of our key observations is the well-known
fact (which can be verified by an easy calculation) that the
spherical derivative is invariant under post-composition with rigid
motions of the Riemann sphere, i.e. that
$$(T\circ f)^\# = f^\#$$
for each $f$ meromorphic in some domain and each rigid motion $T$ of
the sphere. In particular, we have $f^\#=\kl\frac{1}{f}\kr^\#$ for
each meromorphic $f$.

We first prove Theorem \ref{normalresult} for families of {\it
analytic} functions since this is needed in the proof of the
meromorphic case.

\blem\label{normalresult-analytic}
Let some $\varepsilon>0$ be given and set
$$\fa:=\gl f\in\ha(\dk)\;:\; f^\#(z)\ge \varepsilon \mbox{ for all }
z\in\dk\gr.$$
Then $\fa$ is normal in $\dk$.
\elem

We give four different proofs for this fact.

{\bf First Proof of Lemma \ref{normalresult-analytic}.}
First, from $|f'(z)|\ge f^\#(z)\ge \varepsilon$ for all $f\in\fa$ and
all $z\in\dk$ we obtain by Montel's Great Theorem that $\fa'$ is
normal.

Let $(f_n)_n$ be some sequence in $\fa$. After turning to an
appropriate subsequence, we may assume that $(f_n')_n$ converges to
some $d\in\ha(\dk)\cup\gl\infty\gr$.

{\bf Case 1:} $d\in\ha(\dk)$

Then $(f_n')_n$ is locally bounded in $\dk$, and from
\beq\label{estimate}
|f_n^2(z)|\le \frac{|f_n'(z)|}{f_n^\#(z)}\le
\frac{1}{\varepsilon}\cdot |f_n'(z)| \qquad \mbox{ for all } n\in\nat,
z\in\dk
\eeq
we see that $(f_n)_n$ is as well locally bounded in $\dk$, hence
normal there by Montel's theorem.

{\bf Case 2:} $d\equiv \infty$

It suffices to show the normality of $(f_n)_n$ at $z=0$. We fix some
$r\in(0;1)$. W.l.o.g. we may assume that $|f_n'(z)|\ge 1$ for all $z\in
U_r(0)$ and all $n\in\nat$. Then $\log|f_n'|$ is harmonic and positive
in $U_r(0)$ for all $n$. By Harnack's inequality we obtain
$$|f_n'(z)|\le |f_n'(0)|^{\frac{r+|z|}{r-|z|}}\qquad \mbox{ for all }
z\in U_r(0) \mbox{ and all } n\in\nat.$$
Now we choose $\varrho\in(0;r)$ such that
$\frac{r+\varrho}{r-\varrho}<2$. Furthermore, for each $n$ we choose
$z_n\in\dk$ such that $|z_n|=\varrho$ and
$|f_n(z_n)|=\max_{|z|\le\varrho} |f_n(z)|=:M(\varrho,f_n)$. Then from
Cauchy's formula we deduce
$$|f_n'(0)|\le \frac{1}{\varrho}\cdot
M(\varrho,f_n)=\frac{1}{\varrho}\cdot |f_n(z_n)|,$$
hence
$$\varepsilon
\le f_n^\#(z_n)
\le \frac{|f_n'(z_n)|}{|f_n(z_n)|^2}
\le \frac{|f_n'(0)|^{\frac{r+\varrho}{r-\varrho}}}{|f_n(z_n)|^2}
\le \kl\frac{1}{\varrho}\kr^{\frac{r+\varrho}{r-\varrho}}
\cdot |f_n(z_n)|^{\frac{r+\varrho}{r-\varrho}-2}$$
and therefore
$$M(\varrho,f_n)=|f_n(z_n)|
\le\kl\frac{1}{\varepsilon}\cdot\kl\frac{1}{\varrho}
\kr^{\frac{r+\varrho}{r-\varrho}}\kr^{\frac{1}{2-\frac{r+\varrho}{r-\varrho}}}$$
for all $n$. This shows that $(f_n)_n$ is locally uniformly bounded in
$U_\varrho(0)$, hence normal at $z=0$. (And $(f_n')_n$ is uniformly
bounded near $z=0$, too, contradicting $d\equiv\infty$. So, in fact,
Case 2 cannot occur.)

\hfill $\blacksquare$

{\bf Second Proof of Lemma \ref{normalresult-analytic}.}
For each $f\in\fa$ we have $\frac{1}{f}(z)\ne 0$ and
$\l|\kl\frac{1}{f}\kr'(z)\r|
\ge \kl\frac{1}{f}\kr^\#(z)=f^\#(z)\ge\varepsilon,$
hence $\kl\frac{1}{f}\kr'(z)\ne\frac{\varepsilon}{2}$
for all $z\in\dk$. Hence by a well-known normality criterion due to
Y.~Gu \cite{gu}, the family $\gl\frac{1}{f}\;:\;f\in\fa\gr$ is normal
in $\dk$. In view of $f^\#=\kl\frac{1}{f}\kr^\#$ and Marty's
criterion, also $\fa$ is normal in $\dk$. \hfill $\blacksquare$

{\bf Third Proof of Lemma \ref{normalresult-analytic}.}
One can also exploit the Nevanlinna theory\footnote{For the notations
 and main results of Nevanlinna theory, we refer to \cite{hayman-mero}.}
for the proof of Lemma \ref{normalresult-analytic}. As in
(\ref{estimate}), for each $f\in\fa$ and each $z\in\dk$ one has
$$|f^2(z)|\le\frac{1}{\varepsilon}\cdot |f'(z)|.$$
This yields
$$2m(r,f)=m\kl r,f^2\kr
\le m(r,f')+\plog\frac{1}{\varepsilon}
\le m(r,f)+m\kl r,\frac{f'}{f}\kr+\plog\frac{1}{\varepsilon},$$
hence
\beq\label{proof3}
T(r,f)=m(r,f)\le m\kl r,\frac{f'}{f}\kr+\plog\frac{1}{\varepsilon}
\eeq
for each $f\in\fa$ and each $r\in(0;1)$. Now if $\fa$ would be
non-normal in $\dk$, then according to a corollary to the lemma on the
logarithmic derivative due to D.~Drasin \cite[Lemma~6]{drasin}
there would exist a non-normal sequence $(f_n)_n$ in $\fa$, an
$r_0\in(0;1)$ and a constant $C<\infty$ such that
$$m\kl r,\frac{f'_n}{f_n}\kr\le C\cdot\kl \plog
T(r,f_n)+\log\frac{1}{R-r}+1\kr$$
for all $r,R$ satisfying $r_0<r<R<1$ and all $n$. Combining this with
(\ref{proof3}) gives
$$T(r,f_n)\le C\cdot\kl \plog T(r,f_n)+\log\frac{1}{R-r}+1\kr$$
for $r_0<r<R<1$ and all $n$. From this inequality and some standard
arguments (see \cite[p.~118]{schiff}) one easily obtains the normality
of $(f_n)_n$, i.e. a contradiction. Hence $\fa$ is normal in $\dk$. \hfill $\blacksquare$

Our last proof of the analytic case makes us of a famous rescaling
lemma which was originally proved by L.~Zalc\-man~\cite{zalcman} and
later extended by X.-C.~Pang~\cite{pang89,pang90} and by H.~Chen and
Y.~Gu \cite{chengu}.


\blem[Zalcman-Pang Lemma] \label{zalclemma}
Let $\fa$ be a family of meromorphic functions in a domain $D$ all of
whose zeros have multiplicity at least $m$ and all of whose poles have
multiplicity at least $p$ and let $-p<\alpha<m$. Then $\fa$ is not
normal at some $z_0\in D$ if and only if there exist sequences $\kl
f_n\kr_n\In \fa$, $\kl z_n\kr_n \In D$ and $\kl \rho_n\kr_n \In (0,1)$
such that $(\rho_n)_n$ tends to 0, $(z_n)_n$ tends to $z_0$ and such
that the sequence $\kl g_n\kr_n$ defined by
$$g_n(\zeta):= \frac{1}{\rho_n^\alpha}\cdot f_n(z_n+\rho_n \zeta)$$
converges locally uniformly in $\co$ (with respect to the spherical
metric) to a non-constant function $g$ meromorphic in $\co$.
\elem

The case $\alpha=0$ of this lemma (which is the case of the original
Zalcman's Lemma from 1975 \cite{zalcman}) will also be the key tool in
the proof of the meromorphic case of Theorem \ref{normalresult}.

\bbew
The proof of the only-if part can be found in \cite[Lemma
2]{pangzalc2000a}. The converse (which we will use later with
$\alpha=-1$) is probably also known to be true, but,  surprisingly
enough, there  doesn't seem to be an exact reference for it in the
literature (except for the case $\alpha=0$ which is discussed in
\cite{zalcman}). Therefore, we supply the proof.

We assume that $0\le\alpha<m$ and that there exist sequences
$(f_n)_n$, $(z_n)_n$ and $(\rho_n)_n$ as described in the lemma. Then
$\kl g_n^\#\kr_n$ tends to $g^\#$ locally uniformly in $\co$, and
$g^\#\not\equiv0$. 

First we consider the case $0\le \alpha<1$. We take some
$\zeta_0\in\co$ such that $g^\#(\zeta_0)>0$. Then we obtain for all $n$
\beqaro
f_n^\#(z_n+\rho_n\zeta_0)
&=&\frac{|f_n'|}{1+|f_n|^2}(z_n+\rho_n\zeta_0)
=\frac{\rho_n^{\alpha-1}\cdot|g_n'|}{1+\rho_n^{2\alpha}\cdot|g_n|^2}(\zeta_0)\\
&\ge& \rho_n^{\alpha-1}\cdot \frac{|g_n'|}{1+|g_n|^2}(\zeta_0)
=\rho_n^{\alpha-1}\cdot g_n^\#(\zeta_0)\abb \infty \qquad (n\to\infty).
\eeqaro
In view of $\limn (z_n+\rho_n\zeta_0)=z_0$ we conclude by Marty's
criterion that $\fa$ is not normal at $z_0$.

If ($m\ge 2$ and) $1\le \alpha<m$ we set $k:=\lfloor \alpha\rfloor$
and $\beta:=\alpha-k$. Then $0\le \beta<1$. By Weierstrass's theorem
we deduce that $\kl g_n^{(k)}\kr_n$ converges to $g^{(k)}$ locally
uniformly in $\co\mi P_g$. This means that
$\kl\frac{1}{\rho_n^\beta}\cdot f_n^{(k)}(z_n+\rho_n \zeta)\kr_n$
converges to $g^{(k)}$ locally uniformly in $\co\mi P_g$. Here,
$g^{(k)}$ is not constant since otherwise $g$ would be a non-constant
polynomial of degree $\le k<m$, contradicting the fact that all zeros
of $g$ have multiplicity $\ge m$. So by the case $0\le \alpha<1$
already treated (applied with $\beta$ instead of $\alpha$) we conclude
that $\kl f_n^{(k)}\kr_n$ is not normal at $z_0$.

Now we assume that $(f_n)_n$ is normal at $z_0$, say in $U_r(z_0)$ for
some $r>0$. Then after choosing an appropriate subsequence we can
assume that $(f_n)_n$ tends to some $F\in
\ma(U_r(z_0))\cup\gl\infty\gr$ locally uniformly (w.r.t.~the spherical
metric) in $U_r(z_0)$. But $F\equiv \infty$ is impossible since this
would imply $g\equiv\infty$. So $F$ is meromorphic in $U_r(z_0)$, and
hence $\kl f_n^{(k)}\kr_n$ would tend to $F^{(k)}$ locally uniformly
in $U_r(z_0)\mi P_F$. If $F$ would be analytic at $z_0$ we would
obtain a contradiction to the fact that $\kl f_n^{(k)}\kr_n$ is not
normal at $z_0$. Hence $z_0$ has to be a pole of $F$. But then there
exist a $\delta>0$ and an $N\in\nat$ such that $|f_n(z)|\ge 1$ for
each $z\in U_\delta(z_0)$ and each $n\ge N$. By the definition of
$g_n$ this implies $g(\zeta)=\infty$ for each $\zeta\in\co$, a
contradiction. This shows the assertion for the case $\alpha\ge 0$.

If $-m<\alpha\le 0$ we consider the sequence $\kl \frac{1}{g_n}\kr_n$
instead of $(g_n)_n$ and conclude from the case already treated that
$\gl\frac{1}{f}\;|\; f\in\fa\gr$ and hence $\fa$ are not normal.
\ebew

{\bf Fourth Proof of Lemma \ref{normalresult-analytic}.}
Assume that $\fa$ is not normal in $\dk$. Then we apply the
Zalcman-Pang Lemma (Lemma \ref{zalclemma}) with $\alpha=-2$ (which is
admissible since all functions in $\fa$ are analytic) and obtain
sequences $(z_n)_n\In\dk$ and $(\rho_n)_n\In(0;1)$ such that $(z_n)_n$
converges to some $z_0\in\dk$,  $\limn \rho_n=0$ and such that the
sequence $(g_n)_n$ where
$$g_n(\zeta):=\rho_n^2\cdot f_n(z_n+\rho_n\zeta)$$
converges to some nonconstant entire function $g$ locally uniformly in
$\co$. We take some $\zeta_0\in\co$ with $g(\zeta_0)\ne 0$. Then for
all $n$ we have
\beqaro
|g_n'(\zeta_0)|
&=&\l|\rho_n^3 f_n'(z_n+\rho_n\zeta_0)\r|
=\l|\frac{1}{\rho_n}\cdot\frac{f_n'}{f_n^2}(z_n+\rho_n\zeta_0)\cdot g_n^2(\zeta_0)\r|\\
&\ge& \frac{1}{\rho_n}\cdot f_n^\#(z_n+\rho_n\zeta_0)\cdot|g_n^2(\zeta_0)|
\ge \frac{\varepsilon}{\rho_n}\cdot |g_n^2(\zeta_0)|\to\infty \qquad(n\to\infty)
\eeqaro
This contradicts the fact that $g$ is analytic at $\zeta_0$.
Therefore $\fa$ has to be normal. \hfill $\blacksquare$

\section{The Meromorphic Case of Theorem \ref{normalresult}}

We were not able to adjust any of the four approaches presented in
the proof of Lemma \ref{normalresult-analytic} to the meromorphic
case. Our proof of Theorem \ref{normalresult} is based on the
Zalcman-Pang Lemma.

We start with a simple consequence from the argument principle which
might be considered as a counterpart to Weierstrass's theorem.

\blem{}\label{invers-weierstrass}
Let $(f_n)_n$ be some sequence of meromorphic functions in $\dk$ such
that $(f_n')$ converges locally uniformly (w.r.t.~the spherical
metric) to some $d\in\ma(\dk)$. If $(f_n)_n$ converges locally
uniformly (w.r.t.~the spherical or to the euclidean metric) in $\dk\mi
P_d$ to some $F\in \ha(\dk\mi P_d)$, then $F$ is meromorphic in $\dk$
and $(f_n)_n$ converges locally uniformly (w.r.t.~the spherical
metric) in $\dk$ to $F$.
\elem

\bbew{}
Of course, we have $F'=d$ in $\dk\mi P_d$, so the isolated
singularities of $F$ at the points of $P_d$ are poles which means that
$F$ is meromorphic in $\dk$.

First we make sure that $f_n$ is zero-free near the poles of $d$ for
sufficiently large $n$, a fact which is less trivial than it might
seem at first sight.\footnote{In this context it might be instructive
to remember that the locally uniform convergence (w.r.t.~the spherical
metric) cannot be carried over in the other direction, i.e. from
$(f_n)_n$ to $(f_n')_n$; this might fail near the poles of $\limn
f_n$ as the example of the functions
$f_n(z):=\frac{1}{z^2+\frac{1}{n}}$ shows: $(f_n)_n$ converges to
$F(z)=\frac{1}{z^2}$ locally uniformly (w.r.t.~the spherical
metric), but $(f_n'(0))_n$ does not converge to $F'(0)$. In other
words, there is no meromorphic extension of Weierstrass's theorem.}

Let $z_0\in\dk$ be a pole of $d$ and hence of $F$. We choose some
$r>0$ such that $B_r(z_0)\In\dk$, $z_0$ is the only pole of $F$ in
$B_r(z_0)$, $F$ has no zeros in $B_r(z_0)$ and none of the $f_n$ has
zeros or poles on the circle $\partial B_r(z_0)$. By $\pa_h$ resp. $\za_h$
we denote the number of poles resp.~zeros of a meromorphic function
$h$ in $U_r(z_0)$, counting multiplicities. Then by the argument
principle we get
$$\za_{f_n}-\pa_{f_n}
=\frac{1}{2\pi i}\int_{|z-z_0|=r}\frac{f'_n(z)}{f_n(z)}\,dz
\abb \frac{1}{2\pi i}\int_{|z-z_0|=r}\frac{F'(z)}{F(z)}\,dz
=\za_F-\pa_F$$
for $n\to\infty$ since the convergence of $(f_n)_n$ is uniform on
$\partial U_r(z_0)$. So for sufficiently large $n$ we deduce
$$\za_{f_n}-\pa_{f_n}
=\za_F-\pa_F.$$
Here we have $\za_F=0$ by our choice of $r$, and by
Hurwitz's theorem $F'$ and $f_n'$ have the same number of poles in
$U_r(z_0)$ for sufficiently large $n$. Furthermore,
$\pa_{F'}=\pa_F+1$ since $z_0$ is
the only pole of $F$ in $U_r(z_0)$ and trivially $\pa_{f_n'}\ge
\pa_{f_n}+1$ for each $n$. For all $n$ sufficiently large this yields
$$\za_{f_n}=\pa_{f_n}-\pa_F
\le \pa_{f_n'}-1-\kl \pa_{F'}-1\kr=0,$$
hence $\za_{f_n}=0$. So $f_n$ is zero-free in $U_r(z_0)$
for $n$ large enough.

This enables us to apply the maximum principle to
$\kl\frac{1}{f_n}-\frac{1}{F}\kr_n$ in $U_r(z_0)$ to deduce that
$(f_n)_n$ converges to $F$ uniformly (w.r.t. the spherical
metric) in $B_r(z_0)$.

Hence $(f_n)_n$ converges to $F$ locally uniformly (w.r.t.~the
spherical metric) in $\dk$.
\ebew

\bbewo{Theorem \ref{normalresult}}
As in the analytic case, from $|f'(z)|\ge f^\#(z)\ge \varepsilon$ for
all $f\in\fa$ and all $z\in\dk$ we obtain by Montel's Great Theorem
that $\fa'$ is normal. In view of $f^\#=\kl\frac{1}{f}\kr^\#$, we can
also conclude that $\gl\kl\frac{1}{f}\kr'\;\Bigm| \; f\in\fa\gr$ is
normal.

Let $(f_n)_n$ be some sequence in $\fa$. After turning to an
appropriate subsequence, we may assume that $(f_n')_n$ converges to
some $d\in\ma(\dk)\cup\gl\infty\gr$ and that $\kl\frac{1}{f_n}\kr'$
converges to some $\ti{d}\in\ma(\dk)\cup\gl\infty\gr$ locally
uniformly (w.r.t.~the spherical metric) in $\dk$.

{\bf Case 1:} $d\in\ma(\dk)$.

Then $(f_n')_n$ is locally bounded in $\dk\mi P_d$, and from
$$|f_n^2(z)|\le \frac{|f_n'(z)|}{f_n^\#(z)}\le
\frac{1}{\varepsilon}\cdot |f_n'(z)| \qquad \mbox{ for all } n\in\nat,
z\in\dk$$
we see that $(f_n)_n$ is as well locally bounded in $\dk\mi P_d$,
hence normal there by Montel's theorem. W.l.o.g. we may assume that
$(f_n)_n$ converges locally uniformly in $\dk\mi P_d$ to some $F\in
\ha(\dk\mi P_d)$. By Lemma \ref{invers-weierstrass} we conclude that
$(f_n)_n$ converges to $F$ locally uniformly (w.r.t.~the spherical
metric) in $\dk$.

{\bf Case 2:} $\ti{d}\in\ma(\dk)$.

Then the same reasoning as in Case 1 shows that an appropriate
subsequence $\kl\frac{1}{f_{n_k}}\kr_k$ converges to some meromorphic
function locally uniformly in $\dk$. This implies that $(f_{n_k})_k$
converges to some meromorphic function or to $\infty$ locally
uniformly in $\dk$.

{\bf Case 3:} $d\equiv \ti{d}\equiv\infty$.

{\bf Case 3.1:} $\kl\frac{f_n'}{f_n^3}\kr_n$ is normal in $\dk$.

Then w.l.o.g. we may assume that $h_n:=\frac{f_n'}{f_n^3}$ tends to
some $h\in\ma(\dk)\cup\gl\infty\gr$ locally uniformly in $\dk$.

Let some $z_0\in\dk$ be given. Then there exists some $r>0$ such that
$B_r(z_0)\In\dk$ and $h$ omits either the value 0 or the value
$\infty$ in $B_r(z_0)$. So there exists some $N\in\nat$ such that
$h_n$ omits either the value 0 or the value $\infty$ in
$U_r(z_0)$ whenever $n\ge N$. Since the zeros of $f_n$ are poles of
$h_n$ and since the poles of $f_n$ are zeros of $h_n$ we conclude that
$f_n$ omits either the value 0 or the value $\infty$ in $U_r(z_0)$
for $n\ge N$. So for $n\ge N$ either $f_n$ or $\frac{1}{f_n}$ is
analytic in $U_r(z_0)$, and we can apply the analytic version of the
theorem to deduce the normality of $(f_n)_n$ at $z_0$ and hence in the
whole of $\dk$.

{\bf Case 3.2:} $\kl\frac{f_n'}{f_n^3}\kr_n$ is not normal.

Then by the Zalcman-Pang Lemma (Lemma \ref{zalclemma}) with $\alpha=0$
there exist sequences $(z_n)_n\In\dk$ and
$(\rho_n)_n\In(0;1)$ such that $\limn z_n\in\dk$,
$\limn \rho_n=0$ and such that, after replacing $(f_n)_n$ by an
appropriate subseqence, the sequence $(g_n)_n$ where
\beq\label{case32-1}
g_n(\zeta):=\frac{f_n'}{f_n^3}(z_n+\rho_n\zeta)
\eeq
converges to some nonconstant function $g\in\ma(\co)$ locally
uniformly in $\co$.

From our assumption in Case 3 we know that $(f_n')_n$ and
$\kl\frac{f_n'}{f_n^2}\kr_n$ and hence $\kl\frac{f_n'^2}{f_n^2}\kr_n$
tend to $\infty$ locally uniformly in $\dk$. So
$\kl\frac{f_n}{f_n'}\kr_n$ and $\kl\kl\frac{f_n}{f_n'}\kr'\kr_n$ tend
to 0 locally uniformly in $\dk$. Now from
$$\kl\frac{f_n'}{f_n^3}\kr'
=\frac{f_n''}{f_n^3}-3\cdot \frac{f_n'^2}{f_n^4}
=\frac{f_n'^2}{f_n^4}\cdot\kl\frac{f_n f_n''}{f_n'^2}-3\kr
=-\frac{f_n'^2}{f_n^4}\cdot\kl 2+\kl\frac{f_n}{f_n'}\kr'\kr$$
and from
$$g_n'(\zeta):=\rho_n \kl\frac{f_n'}{f_n^3}\kr'(z_n+\rho_n\zeta)\to
g'(\zeta)$$
we see that
\beq\label{case32-2}
-2 \rho_n\cdot \frac{f_n'^2}{f_n^4}(z_n+\rho_n\zeta)\to g'(\zeta)
\qquad (n\to\infty)
\eeq
locally uniformly in $\co\mi P_g$. Since $f_n^\#$ and hence
$\frac{f_n'^2}{f_n^4}$ are non-vanishing, we can apply the maximum
principle to $\kl-\frac{f_n^4}{2\rho_n
  f_n'^2}(z_n+\rho_n\zeta)-\frac{1}{g'(\zeta)}\kr_n$ to deduce that
(\ref{case32-2}) holds also uniformly (w.r.t.~the spherical metric)
near the poles of $g$, hence locally uniformly (w.r.t.~the spherical
metric) in $\co$.

Now from the converse of the Zalcman-Pang Lemma (with
$\alpha=-1$, which is admissible since all poles of
$\frac{f_n'^2}{f_n^4}$ are multiple) we see that $g'$ is constant. So $g$ has the
form $g(\zeta)=A\zeta+B$ with $A\ne0$. Therefore,  choosing an
appropriate branch of the root (and again replacing $(f_n)_n$ by an
appropriate subsequence, if necessary), we have
$$\rho_n^{1/2} \cdot \frac{f_n'}{f_n^2}(z_n+\rho_n\zeta)\to
\sqrt{-\frac{A}{2}} \qquad (n\to\infty),$$
hence
$$\rho_n^{3/2} \cdot\kl\frac{f_n'}{f_n^2}\kr'(z_n+\rho_n\zeta)\to 0 \qquad
(n\to\infty)$$
locally uniformly in $\co$. Inserting this into
$$\kl\frac{f_n'}{f_n^2}\kr'
=\frac{f_n''}{f_n^2}-2\cdot \frac{f_n'^2}{f_n^3}
=\frac{f_n'^2}{f_n^3}\cdot\kl\frac{f_n f_n''}{f_n'^2}-2\kr
=-\frac{f_n'^2}{f_n^3}\cdot\kl 1+\kl\frac{f_n}{f_n'}\kr'\kr,$$
we arrive at
\beq\label{case32-3}
\rho_n^{3/2} \cdot\frac{f_n'^2}{f_n^3}(z_n+\rho_n\zeta)\to 0 \qquad
(n\to\infty)
\eeq
locally uniformly in $\co$. Now, combining (\ref{case32-1}),
(\ref{case32-2}) and (\ref{case32-3}), we obtain
\beqaro
1
&=&\kl\frac{f_n'}{f_n^3}(z_n+\rho_n\zeta)\kr^2
\cdot\kl\rho_n\cdot\frac{f_n'^2}{f_n^4}(z_n+\rho_n\zeta)\kr^{-3}
\cdot\kl\rho_n^{3/2}\cdot\frac{f_n'^2}{f_n^3}(z_n+\rho_n\zeta)\kr^2\\
&&\to g^2(\zeta)\cdot \kl-\frac{g'(\zeta)}{2}\kr^{-3}\cdot 0=0 \qquad (n\to\infty)
\eeqaro
for all $\zeta\in\co$, a contradiction. So this case cannot occur.
\ebewo

\section{Proof of Theorem \ref{universalconstant} and Proposition
  \ref{minimumprinciple}} \label{proofs}

\bbewo{Theorem \ref{universalconstant}}
W.l.o.g. we may assume that $f^\#$ is non-vanishing. (Otherwise,
there's nothing to show.)

There exists a Möbius transformation $T$ which is a rigid motion of the
Riemann sphere such that $T(f(0)) = 0$. We define $g := T\circ f$.
Then $g$ is meromorphic in $\dk$ with $g(0) = 0$, and we have $g^\#(z)
= f^\#(z)$ for all $z\in\dk$ since the spherical derivative is
invariant under rigid motions of the sphere.

The function
$$h(z):=\frac{g(z)}{zg'(z)}$$
is analytic in $\dk$; this follows since the poles of
$g$ are zeros of $h$, $g'$ is non-vanishing and since $g(0) = 0$.
Furthermore, we have $h(0) = 1$.

Let some $r\in(0;1)$ be given. Then the maximum principle implies the
existence of some $z_0 \in\partial U_r(0)$ such that $|h(z_0)|\ge 1$,
and we obtain
$$f^\#(z_0)=g^\#(z_{0})=\frac{| z_{0} g'(z_{0})|}{r\kl1+|g(z_{0})|^2\kr}
=\frac{1}{r\cdot |h(z_0)|}\cdot\frac{| g(z_{0})|}{1+|g(z_{0})|^2} \le
\frac{1}{r}\cdot\frac{1}{2}.$$
This holds for all  $r\in(0;1)$. Therefore, we have
$$\inf_{z\in\dk} f^\#(z)\le \frac{1}{2}.$$
This shows the first assertion.

If $f$ is a rigid motion of the sphere, then so is $g$, and in view of
$g(0)=0$ we deduce that $g(z)=\eta z$ for some $\eta\in\co$ with
$|\eta|=1$. Obviously, $f^\#(z)=g^\#(z)=\frac{1}{1+|z|^2}$ for all
$z\in\dk$, so we have equality in (\ref{infSphericalDeriv}).

Now assume that equality in (\ref{infSphericalDeriv}) holds. Let $T$,
$g$ and $h$ be as in the beginning of the proof. If $h$
would be non-constant, then there would exist some $z^*\in\partial
U_{1/2}(0)$ such that $c:=|h(z^*)|>1$.  Let some $r\in\kl
\frac{1}{2};1\kr$ be given. Then, again by the maximum principle, we
find some $z_0 \in\partial U_r(0)$ such that $|h(z_0)|\ge c$, and as
above we obtain
$$f^\#(z_{0})=\frac{1}{r\cdot |h(z_0)|}\cdot\frac{|g(z_{0})|}{1+|g(z_{0})|^2}
\le \frac{1}{2rc}.$$
Since this holds for all  $r\in\kl\frac{1}{2};1\kr$, we conclude that
$$\inf_{z\in\dk} f^\#(z)\le \frac{1}{2c}<\frac{1}{2},$$
a contradiction to our assumption that $f$ is an extremal function for
(\ref{infSphericalDeriv}). So $h$ must be constant, i.e. $h(z)\equiv
h(0)=1$. But this means that $g(z)=\alpha z$ for some
$\alpha\in\co\mi\gl0\gr$, i.e.
$$\frac{1}{2}
=\inf_{z\in\dk} g^\#(z)
=\inf_{z\in\dk}\frac{|\alpha|}{1+|\alpha z|^2}
=\frac{|\alpha|}{1+|\alpha|^2},$$
and we deduce $|\alpha|=1$. So $g$ and hence $f=T^{-1}\circ g$ are
rigid motions of the sphere.
\ebewo

The use of the function $h$ (which is the crucial idea of the proof)
is inspired by previous work of R.~Fournier and S.~Ruscheweyh
\cite[Theorem~1.2]{fournier-ruscheweyh99} who had used the same method
in the study of certain free boundary value problems.

The proof of Proposition \ref{minimumprinciple} is based on some
elementary facts on conformal metrics and their (Gaussian) curvature
(see \cite[Chapter~2]{krantz}).  Let's recall that the (Gaussian)
curvature $\kappa_\varphi$ of  a metric $\varphi(z)\,|dz|$ is defined
by
$$\kappa_\varphi(z):=-\frac{\Delta (\log\varphi) (z)}{\varphi^2(z)}.$$

\bbewo{Proposition \ref{minimumprinciple}}
Let some $z_0\in\dk$ with $f^\#(z_0)\ne0$ be given. Then, in
particular, $f$ is non-constant. Since $f^\#=\kl\frac{1}{f}\kr^\#$ we
may assume that $z_0$ is no pole of $f$.

So there exists a $\delta>0$ such that $U:=U_\delta(z_0)\In\dk$ and
$f$ is analytic in $U$ and $f^\#$ is zero-free in $U$.

Hence $\lambda(z):=f^\#(z)$ is the density of a conformal metric in
$U$; this metric is the pull-back of the spherical metric
$\frac{|dw|}{1+|w|^2}$ under the analytic map $f|_U$. It is well-known
that the Gaussian curvature of the spherical metric is $+4$ and that
the curvature of a metric coincides with the curvature of its
pull-back under an analytic map. So we conclude that
$$\frac{\Delta (\log f^\#) (z_0)}{\kl f^\#\kr^2(z_0)}=-4,$$
hence
$$\Delta (\log f^\#) (z_0)=-4\cdot\kl f^\#\kr^2(z_0)<0.$$
But this means that $\log f^\#$ cannot have a local minimum at
$z_0$ (since $\Delta (\log f^\#) (z_0)$ is the trace of the Hessian
of $\log f^\#$ at $z_0$, i.e. the sum of its eigenvalues). So $f^\#$
itself cannot have a local minimum at $z_0$.
\ebewo

{\bf Acknowledgment.}
Part of this work was supported by the German Israeli Foundation for
Scientific Research and Development (No.~G~809$-$234.6/2003). Shahar
Nevo also received support from the Israel Science Foundation, Grant
No.~395/07. This research is part of the European Science Foundation
Networking Programme HCAA.

\bibliographystyle{amsplain}

\vspace{12pt}

\vspace{12pt}
\parbox[t]{90mm}{\it
 J\"urgen Grahl \\
 University of W\"urzburg \\
 Department of Mathematics \\
 W\"urzburg\\
 Germany\\
 e-mail: grahl@mathematik.uni-wuerzburg.de}
\parbox[t]{80mm}{\it
Shahar Nevo \\
Bar-Ilan University\\
Department of Mathematics\\
Ramat-Gan 52900\\
Israel\\
e-mail: nevosh@macs.biu.ac.il}

\end{document}